\documentclass[12pt]{amsart}
\usepackage{geometry, amssymb}
\usepackage[all]{xy} 
\title{$K_1$ of some Iwasawa algebras}
\author{Mahesh Kakde}
\date{} 

\newtheorem{theorem}{Theorem}
\newtheorem{lemma}[theorem]{Lemma}
\newtheorem{proposition}[theorem]{Proposition}
\newtheorem{remark}[theorem]{Remark}
\newtheorem{notation}[theorem]{Notations}
\newtheorem{definition}[theorem]{Definition}

\newcommand{\ilim}[1]{%
	\displaystyle{%
	\lim_{\genfrac{}{}{0pt}{}{\longleftarrow}{\scriptstyle #1}} }\;}

\begin{document}

\maketitle
\tableofcontents

\section{Introduction} Let $p$ be a prime. Let $O$ be the ring of integers in a finite extension of $\mathbb{Q}_p$. Let $\mathcal{G}$ be a compact $p$-adic Lie group with a closed normal subgroup $H$ such that $\mathcal{G}/H = \Gamma$ is isomorphic to the additive group of $p$-adic integers. Noncommutative Iwasawa theory has created a of lot interest in the Whitehead group of the Iwasawa algebra $\Lambda_O(\mathcal{G})$ and its localisation at a certain Ore set $S$ (defined by Coates et. al. \cite{CFKSV:2005}. See definition \ref{setS} below). In this paper we describe the Whitehead group of the Iwasawa algebra of $\mathcal{G}$ when the group $\mathcal{G}$ is a pro-$p$ $p$-adic Lie group of dimension one. We also prove some results about the localisation of Iwasawa algebras of such groups. These are generalisations of results in special cases which started with a beautiful work of K. Kato \cite{Kato:2005} (see also K. Kato \cite{Kato:2006},  T Hara \cite{Hara:2009} and M. Kakde\cite{Kakde:2008}). Following K. Kato \cite{Kato:2006} we use the integral logarithm of R. Oliver and L. Taylor as a main tool for proving our results. 

I would like to thank Professor Coates and Professor Burns for their constant encouragement without which I could not have completed this paper.

\subsection{The set up} Let $\mathcal{G}$ be a compact pro-$p$ $p$-adic Lie group of dimension one which has a quotient isomorphic to the additive group of $p$-adic integers. Then $\mathcal{G}$ has a closed normal finite subgroup, say $H$, such that $\mathcal{G}/H = \Gamma \cong \mathbb{Z}_p$. We write $\Gamma$ multiplicatively. Once and for all we choose a lift of $\Gamma$ in $\mathcal{G}$ which gives an isomorphism $\mathcal{G} \cong H \rtimes \Gamma$. Let $\Gamma^{p^e}$ be a fixed open subgroup of $\Gamma$ acting trivially on $H$. Put $G = \mathcal{G}/\Gamma^{p^e}$. Let $O$ be the ring of integers in a finite extension of $\mathbb{Q}_p$. The Iwasawa algebra of $\mathcal{G}$ with coefficients in $O$ is defined as
\[
\Lambda_O(\mathcal{G}) = \ilim{U} O[\mathcal{G}/U],
\]
where the inverse limit is over all open normal subgroups $U$ of $\mathcal{G}$. 

\subsection{Twisted group rings} Recall the definition of twisted group rings. Let $R$ be a ring and $P$ be any group. Let 
\[
\tau: P \times P \rightarrow R^{\times},
\]
be a two cocycle. Then the twisted group ring, denoted by $R[P]^{\tau}$, is a free $R$-module generated by $P$. We denote the image of $p \in P$ in $R[P]^{\tau}$ by $\overline{p}$. Hence, every element of $R[P]^{\tau}$ is a finite sum $\sum_{p \in P} r_p \overline{p}$ and the addition is component wise. The multiplication has the following twist: 
\[
\overline{p}\cdot \overline{p'} = \tau(p,p') \overline{pp'}.
\]

\subsection{The Iwasawa algebra as a twisted group ring} The Iwasawa algebra $\Lambda_{O}(\mathcal{G})$ is a twisted group ring
\[
\Lambda_{O}(\mathcal{G}) = \Lambda_{O}(\Gamma^{p^e})[G]^{\tau},
\]
where $\tau$ is the twisting map given by 
\[
\tau(h_1\gamma^{a_1}, h_2\gamma^{a_2}) = \gamma^{[\frac{a_1+a_2}{p^e}]p^e} \in \Lambda_{O}(\Gamma^{p^e})^{\times},
\]
where $[x]$ denotes the greatest integer less than or equal to $x$. Note that for any $g, g' \in G$, 
\[
\tau(g, g') = \tau(g',g) \qquad \text{i.e. $\tau$ is \emph{symmetric}},
\]
\[
\tau(g,g^{-1})=1=\tau(g,1).
\]

\subsection{Definition of the Ore set $S$} \label{setS} Recall the Ore set $S$, due to Coates et. al. \cite{CFKSV:2005}, mentioned in the introduction:
\[
S = \{ f \in \Lambda_O(\mathcal{G}) | \Lambda_O(\mathcal{G})/\Lambda_O(\mathcal{G})f \text{ is a finitely generated } O-\text{module}\}
\]
It is proven in \emph{loc. cit.} that $S$ is a multiplicatively closed left and right Ore set and does not contain any zero divisors. Hence we may localise to get the ring $\Lambda_O(\mathcal{G})_S$ which contains $\Lambda_O(\mathcal{G})$.

\subsection{Some useful lemmas}

\begin{lemma} The subset of $S$ defined by $T=\Lambda_O(\Gamma^{p^e}) - p\Lambda_O(\Gamma^{p^e})$ is a multiplicatively closed left and right Ore subset of $\Lambda_O(\mathcal{G})$. The inclusion of rings $\Lambda_O(\mathcal{G})_T \rightarrow \Lambda_O(\mathcal{G})_S$ is an isomorphism.
\label{lemma1}
\end{lemma}
\noindent{\bf Proof:} Since the group $\Gamma^{p^e}$ is central in $G$, it is clear that $T$ is a left and right Ore set. Since $\Lambda_O(\Gamma^{p^e})$ is a domain, the set $T$ does not contain any zero divisors. Hence, the map $\Lambda_O(\mathcal{G})_T \rightarrow \Lambda_O(\mathcal{G})_S$, induced by the inclusion $T \rightarrow S$, is an inclusion. We now prove that it is surjective.

Note that $\Lambda_O(\mathcal{G})_{T} = \Lambda_O(\Gamma^{p^e})_T \otimes_{\Lambda_O(\Gamma^{p^e})} \Lambda_O(\mathcal{G})$. We first show that 
\[
Q(\Lambda_O(\mathcal{G})) = Q(\Lambda_O(\Gamma^{p^e})) \otimes_{\Lambda_O(\Gamma^{p^e})} \Lambda_O(\mathcal{G}),
\] 
where $Q(R)$ denotes the total ring of fractions of a ring $R$. Note that we have an injective map
\[
Q(\Lambda_O(\Gamma^{p^e})) \otimes_{\Lambda_O(\Gamma^{p^e})} \Lambda_O(\mathcal{G}) \hookrightarrow Q(\Lambda_O(\mathcal{G})).
\]
As $Q(\Lambda_O(\Gamma^{p^e}))$ is a field and $\Lambda_O(\mathcal{G})$ is a free $\Lambda_O(\Gamma^{p^e})$-module of finite rank, $Q(\Lambda_O(\Gamma^{p^e})) \otimes_{\Lambda_O(\Gamma^{p^e})} \Lambda_O(\mathcal{G})$ is an Artinian ring. Hence every regular element is invertible. The ring $\Lambda_O(\mathcal{G})$ is contained in $Q(\Lambda_O(\Gamma^{p^e})) \otimes_{\Lambda_O(\Gamma^{p^e})} \Lambda_O(\mathcal{G})$ and every regular element of $\Lambda_O(\mathcal{G})$ is invertible in $Q(\Lambda_O(\Gamma^{p^e})) \otimes_{\Lambda_O(\Gamma^{p^e})} \Lambda_O(\mathcal{G})$, hence the injection $Q(\Lambda_O(\Gamma^{p^e})) \otimes_{\Lambda_O(\Gamma^{p^e})} \Lambda_O(\mathcal{G}) \hookrightarrow Q(\Lambda_O(\mathcal{G}))$ must be surjective. Any element $x \in \Lambda_O(\mathcal{G})_S \subset Q(\Lambda_O(\mathcal{G}))$ can be written as $\frac{a}{t}$, with $a \in \Lambda_O(\mathcal{G})$ and $ t $ a non-zero element of $\Lambda_O(\Gamma^{p^e})$. If $t \in p^n\Lambda_O(\Gamma^{p^e})$, then $tx = a \in p^n\Lambda_O(\mathcal{G})_S$. On the other hand, $a$ also lies in $\Lambda_O(\mathcal{G})$. Hence $a \in p^n\Lambda_O(\mathcal{G})$. Then we can divide the largest possible power of $p$ from $t$ and the same power of $p$ from $a$, and $x$ can be presented as $\frac{a}{t}$, with $a \in \Lambda_O(\mathcal{G})$ and $ t \in T$. \qed

\begin{remark} We remark that $\Lambda_O(\mathcal{G})_T = \Lambda_O(\Gamma^{p^e})_T[G]^{\tau}$, for the same twisting map $\tau$ as above. We also study the $p$-adic completion $\widehat{\Lambda_{O}(\mathcal{G})_{T}} = \widehat{\Lambda_O(\Gamma^{p^e})_T}[G]^{\tau}$ of $\Lambda_O(\mathcal{G})_T$. Note that $\Lambda_O(\Gamma^{p^e})_T = \Lambda_O(\Gamma^{p^e})_{(p)}$.
\end{remark}

\begin{notation} In the rest of this section $R$ denotes either $\Lambda_O(\Gamma^{p^e})$ or $\widehat{\Lambda_O(\Gamma^{p^e})_{(p)}}$. For any subgroup $P$ of $G$, we denote the inverse image of $P$ in $\mathcal{G}$ by $U_P$. Let $N_GP$ be the normaliser of $P$ in $G$. We put $W_GP= N_GP/P$, the Weyl group. We denote the set of all cyclic subgroups of $G$ by $C(G)$.
\end{notation}

\begin{definition} Let $I_1$ and $I_2$ be ideals in $R[G]^{\tau}$. Then we define $[I_1, I_2]$ to be the additive subgroup of $R[G]^{\tau}$ generated by elements of the form $ab-ba$, where $a \in I_1$ and $b \in I_2$.
\end{definition}

\begin{lemma} We have the isomorphism of $R$-modules
\[
R[G]^{\tau}/[R[G]^{\tau}, R[G]^{\tau}] \rightarrow R[Conj(G)]^{\tau}.
\]
\end{lemma}

\begin{remark} Throughout this paper we usually denote the class of $g$ in $Conj(G)$ by $[g]$.
\end{remark}

\noindent{\bf Proof:} Consider the $R$-module homomorphism 
\[
R[G]^{\tau} \xrightarrow{\phi} R[Conj(G)]^{\tau},
\]
\[
\sum r_g \overline{g} \mapsto \sum r_g [\overline{g}].
\]
This map $\phi$ is surjective and since $\tau$ is symmetric, the kernel of $\phi$ contains $[R[G]^{\tau}, R[G]^{\tau}]$. We must show that it is equal to $[R[G]^{\tau}, R[G]^{\tau}]$. Let $\sum r_g\overline{g} \in ker(\phi)$. Let $C_g$ denote the centraliser of $g$ in $G$.  Then for each $g \in G$, we get
\[
\sum_{x \in G/C_g} r_{xgx^{-1}} = 0.
\]
Consider
\begin{align*}
\sum r_{xgx^{-1}} \overline{xgx^{-1}} - \sum r_{xgx^{-1}} \overline{g} & = \sum r_{xgx^{-1}}(\overline{xgx^{-1}} - \overline{g}) \\
& = \sum r_{xgx^{-1}} (\tau(xg, x^{-1})^{-1} \overline{xg} \overline{x^{-1}} - \tau(x^{-1}, xg)^{-1} \overline{x^{-1}}\overline{xg}) \\
& = \sum r_{xgx^{-1}} \tau(xg, x^{-1})^{-1} ( \overline{xg}\overline{x^{-1}} - \overline{x^{-1}} \overline{xg}),
\end{align*}
which clearly belongs to $[R[G]^{\tau}, R[G]^{\tau}]$. All the sums are over $x \in G/C_g$. \qed

\begin{remark} Here of course $R[Conj(G)]^{\tau}$ is just a $R$-module and multiplication is not defined. However, taking powers is well defined and while doing so we must remember the twisting map $\tau$. Hence we put it in the notation.
\end{remark}

\begin{lemma} For any subgroup $P \leq G$, we have
\[
\Lambda_O(U_P) \cong \Lambda_O(\Gamma^{p^e})[P]^{\tau},
\]
and
\[
\widehat{\Lambda_O(U_P)_S} \cong \widehat{\Lambda_O(\Gamma^{p^e})_{(p)}}[P]^{\tau},
\]
for the same twisting map $\tau$ as above.
\end{lemma}

\begin{lemma} If $P \leq G$ is a cyclic subgroup, then $U_P$ is abelian (though $U_P$ is not necessarily a direct product of $\Gamma^{p^e}$ and $P$). 
\end{lemma}

\section{The additive side} 

The group $N_GP$ acts $R$-linearly on $R[P]^{\tau}$ by conjugation on $P$. Let $C(G)$ be the set of cyclic subgroups of $G$. For any $P \in C(G)$, define a map 
\[
t^G_P: R[Cong(G)]^{\tau} \rightarrow R[P]^{\tau}
\]
as follows: let $C(G,P)$ denote any set of left coset representatives of $P$ in $G$. Then 
\[
t^G_P(\overline{g}) = \sum_{x \in C(G,P)} \{ \overline{x}_i^{-1}\overline{g}\overline{x}_i |  x_i^{-1}gx_i \in P\}.
\]
This is a well-defined $R$-linear map, independent of the choice of $C(G,P)$. For any $P \in C(G)$, define
\[
\eta_P : R[P]^{\tau} \rightarrow R[P]^{\tau},
\]
by 
\[
\eta_P(h)=  \left\{ 
\begin{array}{l l}
h & \quad \mbox{if $h$ is a generator of $P$} \\
0 & \quad \mbox{if not}. \\
\end{array} \right.
\]

\begin{definition} Define $\beta^G_P : R[Conj(G)]^{\tau} \rightarrow R[P]^{\tau}$ by $\beta^G_P = \eta_P \circ t^G_P$ and $\beta^G_R$ by 
\[
\beta^G_R = (\beta^G_P)_{P\in C(G)} : R[Conj(G)]^{\tau} \rightarrow \prod_{P\in C(G)}R[P]^{\tau}.
\]
Sometimes we denote $\beta^G_{\Lambda_O(\Gamma^{p^e})}$ by just $\beta^G$.
\end{definition}

\begin{definition} Let $P \in C(G)$ be a cyclic subgroup of $G$. We define $T_{P,R}$ to be the image of the map 
\[
tr: R[P]^{\tau} \rightarrow R[P]^{\tau} \qquad x \mapsto \sum_{g \in W_GP} \overline{g} \overline{x} \overline{g}^{-1}.
\]
It is an ideal in the ring $(R[P]^{\tau})^{W_GP}$. We denote $T_{P, \Lambda_O(\Gamma^{p^e})}$ by $T_P$ and $T_{P, \widehat{\Lambda_O(\Gamma^{p^e})_{(p)}}}$ by $\widehat{T_P}$.
\end{definition}

\begin{definition} Let $P \leq P_1$ be two cyclic subgroups of $G$. Then we have the norm map
\[
nr^{P_1}_P : (R[P_1]^{\tau})^{\times} \rightarrow (R[P]^{\tau})^{\times},
\]
and the trace map
\[
tr^{P_1}_P : R[P_1]^{\tau} \rightarrow R[P]^{\tau}.
\]
\end{definition}

\begin{definition} Let $\psi^G_R \subset \prod_{P \in C(G)} R[P]^{\tau}$ be the subgroup consisting of all tuples $(a_P)$ such that \\
A1. If $P < P_1 \in C(G)$, then $tr^{P_1}_P(a_{P_1}) = 0$. \\
A2. $(a_P)_{P\in C(G)}$ is invariant under conjugation action by every $g \in G$. \\
A3. For all $P \in C(G)$, $a_P \in T_{P,R}$. \\
Sometimes we denote $\psi^G_{\Lambda_O(\Gamma^{p^e})}$ by just $\psi^G$.
\end{definition}

\subsection{The additive theorem}

\begin{theorem} The homomorphism $\beta^G_R$ induces an isomorphism between $R[Conj(G)]^{\tau}$ and $\psi^G_R$.
\label{additivetheorem}
\end{theorem}

\begin{lemma} The image of $\beta^G_R$ is contained in $\psi^G_R$. 
\label{additivecontain}
\end{lemma}
\noindent{\bf Proof:} It is enough to show that $\beta^G_R([\overline{g}]) \in \psi^G_R$, for any $g \in G$, i.e. it satisfies A1, A2 and A3. \\

\noindent A1. For $P < P_1 \in C(G)$, it is clear that $tr^{P_1}_P([\bar{g}]) = 0$ unless $ g \in P$, in which case $tr^{P_1}_P([\bar{g}])  = [P_1:P][\bar{g}]$. Since the coefficient of any element $g' \in P$ in $\beta^G_{P_1}([\overline{g}])$ is 0 if $g'$ does not generate $P_1$, it is clear that $tr^{P_1}_P(\beta^G_{P_1}([\overline{g}])) = 0$.  \\
 
\noindent A2. For any $ g, g_1 \in G$, we must show that $\overline{g_1} \beta^G_P([\overline{g}]) \overline{g_1^{-1}} = \beta^G_{g_1Pg_1^{-1}}([\overline{g}])$, for any $P \in C(G)$. 

\begin{align*} 
\overline{g_1}t^G_P([\overline{g}])\overline{g_1^{-1}} & = \overline{g_1}(\sum_{x \in C(G,P)} \{ [\overline{x^{-1}gx}] : x^{-1}gx \in P\}) \overline{g_1^{-1}} \\
& = \sum_{x \in C(G,P)} \{ [\overline{ (g_1x^{-1}g_1^{-1})(g_1gg_1^{-1})(g_1xg_1^{-1})}] : x^{-1}gx \in P\} \\
& = \sum_{x_1 \in C(G, g_1Pg_1^{-1})} \{[\overline{x_1^{-1}(g_1gg_1^{-1})x_1}] : x_1^{-1}g_1gg_1^{-1}x_1 \in g_1Pg_1^{-1}\} \\
& = t^G_{g_1Pg_1^{-1}}([\overline{g_1gg_1^{-1}}]) \\
& = t^G_{g_1Pg_1^{-1}}([\overline{g}]).
\end{align*}
Note that above we have, on several occasions, used $\bar{g_1} \bar{g} \overline{g_1^{-1}} = \overline{g_1gg_1^{-1}}$. Hence $g_1\eta_P(t^G_P([\overline{g}]))g_1^{-1} = \eta_{g_1Pg_1^{-1}}(t^G_{g_1Pg_1^{-1}}([\overline{g}]))$. \\

\noindent A3. We must show that $t^G_P([\bar{g}]) \in T_{P,R}$ for any $P \in C(G)$. 
\[
t^G_P([\overline{g}]) = t^{N_GP}_P(t^G_{N_GP}([\overline{g}])).
\]
Hence, it is enough to show that $t^{N_GP}_P([\overline{g}]) \in T_{P,R}$ for any $g \in N_GP$. But $t^{N_GP}_P([\overline{g}])$ is non-zero if and only if $ g \in P$, and when $ g \in P$, we have
\begin{align*}
t^{N_GP}_P([\overline{g}]) & = \sum_{x \in C(N_GP, P)} \overline{x^{-1}gx} \\
& = \sum_{x \in W_GP} \overline{x^{-1}gx} \in T_{P,R}.
\end{align*}
This finishes proof of the lemma. \qed

\begin{definition} We define a left inverse $\delta$  of $\beta^G_R$ by
\[
\delta :  \prod_{P\in C(G)} R[P]^{\tau} \rightarrow R[Conj(G)]^{\tau}[\frac{1}{p}],
\]
by putting $\delta = \sum_{P \in C(G)} \delta_P$ and defining $\delta_P$ by
\[
\delta_P : R[P]^{\tau} \rightarrow R[Conj(G)]^{\tau}[\frac{1}{p}],
\]
\[
x \in R[P]^{\tau} \mapsto \frac{1}{[G:P]}[x] \in R[Conj(G)]^{\tau}[\frac{1}{p}].
\]
\end{definition}

\begin{lemma} $\delta \circ \beta^G_R$ is identity on $R[Conj(G)]^{\tau}$. In particular, $\beta^G_R$ is injective.
\end{lemma}

\noindent{\bf Proof:} For any $g \in G$, we show that $\delta(\beta^G_R([\overline{g}])) = [\overline{g}]$. Let $P$ be the cyclic subgroup of $G$ generated by $g$. Let $C$ be the set of all conjugates of $P$ in $G$. Then
\begin{align*} 
\delta(\beta^G_R([\overline{g}])) & = \sum_{P' \in C} \delta_{P'}(\beta^G_R([\overline{g}])) \\
&= \sum_{P'\in C} \frac{1}{[G:P]} [\beta^G_R([\overline{g}])]\\
&= \frac{1}{[G:P]}\sum_{P' \in C} [N_GP':P'] [\overline{g}] \\
&= \frac{1}{[G:N_GP]} \sum_{P' \in C} [\overline{g}] = [\overline{g}]
\end{align*}

\qed

\begin{lemma} The restriction of $\delta$ to the subgroup $\psi^G_R$ is injective and its image lies in $R[Conj(G)]^{\tau}$.
\end{lemma} 

\noindent{\bf Proof:} Let $(a_P) \in \psi^G_R$ be such that $\delta((a_P)) = 0$. We claim that $\delta_P(a_P) = 0$ for each $P \in C(G)$. This follows from two simple observations: firstly, by A1 $\delta_P(a_P)$ and $\delta_{P'}(a_{P'})$ cannot cancel each other unless $P$ and $P'$ are conjugates; but when $P$ and $P'$ are conjugates, $\delta_P(a_P) = \delta_{P'}(a_{P'})$ by A2. Hence $\delta_P(a_P) =0$ for every $P \in C(G)$. 

Let $a_P = \sum_{g \in P} r_g \overline{g}$. Then $\delta_P(a_P) = \frac{1}{[G:P]} \sum_{g \in P} r_g [\overline{g}]$. Let $H_0(N_GP,P)$ be the orbit set for the conjugate action of $N_GP$ on $P$. Then
\[
\delta_P(a_P) = \sum_{x \in H_0(N_GP,P)} \sum_{g \in x} r_g[\overline{g}].
\]
If $g, g' \in x \in H_0(N_GP,P)$, then $r_g = r_{g'}$ by A2. Call it $r_x$. Hence, 
\[
\delta_P(a_P) = \sum_{x \in H_0(N_GP, P)} r_x \sum_{g \in x} [\overline{g}] = 0.
\]
Hence $r_x = 0$ for all $x \in H_0(N_GP, P)$.

A3 says that $a_P \in T_{P,R}$ for every $P \in C(G)$. Let $a_P = tr(b_P)$ for some $b_P \in R[P]^{\tau}$. Then
\begin{align*}
\delta_P(a_P) & = \delta_P(\sum_{x \in W_GP} xb_Px^{-1}) \\
& = [N_GP:P] \delta_P(b_P).
\end{align*}
On the other hand
\begin{align*}
\sum_{x \in C(G, N_GP)} \delta_{xPx^{-1}}(a_{xPx^{-1}}) & = [G:N_GP] \delta_P(a_P), \qquad \text{by A2} \\
& = [G:P] \delta_P(b_P) \in R[Conj(G)]^{\tau}.
\end{align*}

\qed
 
\noindent{\bf Proof of theorem \ref{additivetheorem}:} $\delta|_{\psi^G_R}$ is injective and $\delta \circ \beta^G_R$ is identity on $R[Conj(G)]^{\tau}$. We claim that $\beta^G_R \circ \delta$ is identity on $\psi^G_R$. Let $(a_P) \in \psi^G_R$, then $\delta(\beta^G_R(\delta((a_P)))) = \delta((a_P))$. Since the image of $\beta^G_R$ is contained in $\psi^G_R$ and $\delta$ is injective on $\psi^G_R$, we get $\beta^G_R(\delta((a_P))) = (a_P)$. 

\qed

\section{Definition of $K_1$} In this section we recall the definition of $K_1$ groups and take it as an opportunity to introduce notations. Let $R$ be any associative ring with a unit. For each integer $n > 0$, let $GL_n(R)$ denote the group of invertible $n \times n$ matrices with entries in $R$. Regard $GL_n(R)$ as a subgroup of $GL_{n+1}(R)$ by identifying $A \in GL_n(R)$ with $\left( \begin{array}{l l} A & 0 \\ 0 & 1 \end{array} \right) \in GL_{n+1}(R)$. Set $GL(R) = \cup_{n=1}^{\infty} GL_n(R)$. For any $n$, any $i \neq j$, with $1 \leq i,j \leq n$, and any $r \in R$, let $e^r_{ij} \in GL_n(R)$ be an elementary matirx i.e. which has 1's on the diagonal and $r$ in the $(i,j)th$ position. For each integer $n > 0$, let $E_n(R) \subset GL_n(R)$ be the subgroup generated by all elementary $n\times n$ matrices. Set $E(R) = \cup_{n=1}^{\infty} E_n(R)$. Whitehead lemma says that $E(R) = [GL(R), GL(R)]$, the commutator subgroup of $GL(R)$. In particular, $E(R)$ is normal in $GL(R)$ and the quotient $GL(R)/E(R)$ is abelian which we denote by $K_1(R)$.
 
Let $I \subset R$ be any ideal. Denote the group of invertible matrices which are congruent to the identity modulo $I$ by $GL(R,I)$ and denote the smallest normal subgroup of $G(R)$ containing all $e_{ij}^r$ for all $r \in I$ by $E(R,I)$. Finally, set $K_1(R,I) = GL(R,I)/E(R,I)$. The lemma of Whitehead says that $E(R,I) = [GL(R), GL(R,I)]$. Hence $K_1(R,I)$ is an abelian group.

\section{Logarithm and integral logarithm for Iwasawa algebras} 

\subsection{Logarithm homomorphism} Logarithms on $K_1$ groups of $p$-adic orders were constructed by R. Oliver and L. Taylor. In this section we follow R. Oliver \cite{Oliver:1988} to define logarithm homomorphism on $K_1$-groups of Iwasawa algebras, a straightforward generalisation of the construction of R. Oliver and L. Taylor. As before let $R$ denote the ring $\Lambda_O(\Gamma^{p^e})$ or $\widehat{\Lambda_O(\Gamma^{p^e})_{(p)}}$. Let $J_R$ be the Jacobson radical of $R[G]^{\tau}$. Since $\mathcal{G}$ is pro-$p$, the ring $R[G]^{\tau}$ is a local ring and hence $J_R$ is its maximal ideal. We have the series
\[
Log(1+x) = \sum_{i=1}^{\infty} (-1)^{i+1}\frac{x^i}{i},
\]
and 
\[
Exp(x) = \sum_{i=0}^{\infty} \frac{x^i}{i!}.
\]

\begin{lemma} $J_R/pR[G]^{\tau}$ is a nilpotent ideal of $R[G]^{\tau}/pR[G]^{\tau}$.
\end{lemma}

\noindent{\bf Proof:} Let $k = O/(p)$. We have the following exact sequence
\[
0 \rightarrow J_R/pR[G]^{\tau} \rightarrow Q(k[[\mathcal{G}]]) \rightarrow Q(k[[\Gamma]]) \rightarrow 0.
\]
Let $N$ be the kernel of the map $k[[\mathcal{G}]] \rightarrow k[[\Gamma]]$, and let $I_H$ be the kernel of the map $k[H] \rightarrow k$. Then $N = k[[\mathcal{G}]]I_H$, and since $H$ is a finite $p$ group, we have
\[
N^n = k[[\mathcal{G}]] I_H^n = 0,
\]
for some positive integer $n$. By lemma \ref{lemma1} we can write any element $x \in Q(k[[\mathcal{G}]])$ as $x = \frac{a}{t}$, with $a \in k[[\mathcal{G}]]$ and $t \in k[[\Gamma^{p^e}]]$. Also, $x \in J_R/pR[G]^{\tau}$ if and only if $a \in N$. As $t$ is central, we deduce that $J_R/pR[G]^{\tau}$ is nilpotent. 

\qed

\begin{lemma} Let $I \subset J_R$ be any ideal of $R[G]^{\tau}$. Then \\ 
1) For any $x\in I$, the series $Log(1+x)$ converges to an element in $R[G]^{\tau}[\frac{1}{p}]$. Moreover, for any $u, v \in 1+I$
\begin{equation}
\label{log1}
Log(uv) \equiv Log(u) + Log(v) (\text{mod } [R[G]^{\tau}[\frac{1}{p}], I[\frac{1}{p}]]).
\end{equation}
2) If $I \subset \xi R[G]^{\tau}$, for some central element $\xi$ such that $\xi^p \in p\xi R[G]^{\tau}$, then for all $u,v \in 1+I$, $Log(u)$ and $Log(v)$ converge to an element in $I$ and
\begin{equation}
\label{log2}
Log(uv) \equiv Log(u) + Log(v) (\text{mod } [R[G]^{\tau}, I]).
\end{equation}
In addition, if $I^p \subset pIJ_R$, then the series $Exp(x)$ converges to an element in $1+I$ for all $x \in I$; the maps defined by $Exp$ and $Log$ are inverse bijections between $I$ and $1+I$. Moreover, $Exp([R[G]^{\tau}, I]) \subset E(R[G]^{\tau},I)$ and for any $x,y \in I$, we have
\begin{equation}
\label{exp}
Exp(x+y) \equiv Exp(x)\cdot Exp(y) (\text{mod } E(R[G]^{\tau},I)).
\end{equation}
 
\end{lemma}

\noindent{\bf Proof:} This is analogue of lemma 2.7 in R. Oliver \cite{Oliver:1988}. We reproduce the proof here stressing the modifications needed in our situation. \\

\noindent{\bf Step 1.} For any $n \geq 1$, the ideal $J_R/p^nR[G]^{\tau}$ is nilpotent in $R[G]^{\tau}/p^nR[G]^{\tau}$. Hence for any $x \in I \subset J_R$, the terms $x^n/n$ converge to 0 and the series $Log(1+x)$ converges in $I[\frac{1}{p}] \subset R[G]^{\tau}[\frac{1}{p}]$.

If $I \subset \xi R[G]^{\tau}$, for some central $\xi$ such that $\xi^p \in p\xi R[G]^{\tau}$, then $I^p \subset pI$. Hence $I^n \subset nI$ for every positive integer $n$. So $x^n/n \in I$. Hence the series $Log(1+x)$ converges to an element in $I \subset R[G]^{\tau}$.

Furthermore, if $I^p \subset pIJ_R$, note that for any positive integer $n$ such that $p^k \leq  n < p^{k+1}$
\[
I^n \subset p^{([n/p]+[n/p^2]+\cdots+[n/p^k])}I J_R^k = n! IJ^k_R.
\]
Recall that $[y]$ denotes the greatest integer less than or equal to $y$ and that $n! p^{-([n/p]+\cdots+[n/p^k])}$ is a $p$-adic unit. Hence $Exp(x)$ converges to an element in $1+I$. The fact that $Log$ and $Exp$ are inverse bijections between $1+I$ and $I$ is formal. \\

\noindent{\bf Step 2.} For any $I \subset J_R$, set 
\[
U(I) = \sum_{m \geq 0, n \geq 1} \frac{1}{m+n} [I^m,I^n] \subset [R[G]^{\tau}[\frac{1}{p}], I[\frac{1}{p}]],
\]
a $R$-submodule of $R[G]^{\tau}[\frac{1}{p}]$. If $I \subset \xi R[G]^{\tau}$, where $\xi$ is a central element such that $\xi^p \in p\xi R[G]^{\tau}$, then $\xi^n \in n\xi R[G]^{\tau}$, and 
\[
U(I) = \langle [r, \frac{\xi^{m+n}}{m+n}s] : m \geq 0, n \geq 1, \xi^mr \in I^m, \xi^ns\in I^n, \xi r, \xi s \in I \rangle \subset [R[G]^{\tau}, I].
\]
So the congruences (\ref{log1}) and (\ref{log2}) will both follow once we have shown that for every $I \subset J_R$ and every $x,y \in I$
\[
Log((1+x)(1+y)) \equiv Log(1+x) + Log(1+y) (\text{mod } U(I)).
\]

For each $n \geq 1$, we let $W_n$ be the set of formal ordered monomials of length $n$ in two variables $a,b$. For $w \in W_n$, set 

$C(w)$ = orbit of $w$ in $W_n$ under cyclic permutations. 

$k(w)$ = number of occurrences of $ab$ in $w$.

$r(w)$ = coefficients of $w$ in $Log(1+a+b+ab) = \sum_{i=0}^{k(w)} (-1)^{n-i-1} \frac{1}{n-i}\binom{k(w)}{i}$.

If $w^{\prime} \in C(w)$, then it is clear that $w(x,y) \equiv w^{\prime}(x,y) (mod \ [I^i, I^j])$ for some $i,j \geq 1$ such that $i+j=n$. So
\begin{align*}
Log(1+x+y+xy) & = \sum_{n=1}^{\infty} \sum_{w\in W_n} r(w)w(x,y) \\
& \equiv \sum_{n=1}^{\infty} \sum_{w \in W_n/C} \Big( \sum_{w^{\prime} \in C(w)} r(w^{\prime})\Big)w(x,y) (mod \ U(I)).
\end{align*}
Let $ k =max \{k(w^{\prime}) : w^{\prime} \in C(w)\}$. Let $ |C(w)| = n/t$. Then $C(w)$ contains $k/t$ elements with exactly $(k-1)$ $ab$'s and $(n-k)/t$ elements with $k$ $ab$'s. Hence
\begin{align*}
\sum_{w^{\prime} \in C(w)} r(w^{\prime}) &= \frac{1}{t} \sum_{i=0}^{k} (-1)^{n-i-1} \frac{1}{n-i}\Big((n-k)\binom{k}{i} + k \binom{k-1}{i}\Big) \\
                                                                        &=\frac{1}{t} \sum_{i=0}^{k} (-1)^{n-i-1} \frac{1}{n-i}\Big((n-k)\binom{k}{i} + (k-i)\binom{k}{i}\Big) \\
                                                                        &= \frac{1}{t}\sum_{i=0}^{k} (-1)^{n-i-1} \binom{k}{i},
\end{align*}
which is 0 unless $k=0$, in which case it is equal to $(-1)^{n-1}\frac{1}{n}$. Thus
\[
Log(1+x+y+xy) \equiv \sum_{n=1}^{\infty} (-1)^{n-1} \Big(\frac{x^n}{n}+\frac{y^n}{n}\Big) = Log(1+x)+Log(1+y) \ (mod\ U(I)).
\]

\noindent{\bf Step 3.} We now prove the congruence (\ref{exp}). $Exp$ and $Log$ induce bijection between $I$ and $1+I$. Hence 
\[
Log(Exp(x)Exp(y)) \equiv x+y \ (mod \ U(I)),
\]
which gives 
\begin{align*}
Exp(x)Exp(y)Exp(x+y)^{-1} & \in Exp(x+y+U(I))Exp(-x-y) \\
& \subset Exp(U(I)) \subset Exp([R[G]^{\tau},I]).
\end{align*}

Hence we only need to prove that $Exp([R[G]^{\tau},I])$ is contained in $E(R[G]^{\tau},I)$. Choose a $R$-basis $\{[s_1,v_1], \ldots, [s_m,v_m]\}$ of $[R[G]^{\tau}, I]$, with $s_i \in R[G]^{\tau}$ and $v_i \in I$. Let $x = \sum_{i=1}^{m}a_i[s_i,v_i]$ be an element in $[R[G]^{\tau},I]$. Define 
\[
\psi(x)= \prod_{i=1}^{m}(Exp(a_is_iv_i)Exp(a_iv_is_i)^{-1}).
\]
For any $r \in R[G]^{\tau}$ and any $x \in I$, an identity of Vaserstein gives
\begin{align*}
Exp(rx)Exp(xr)^{-1} & = \Big( 1+ r(\sum_{n=1}^{\infty} \frac{x(rx)^{n-1}}{n!})\Big)\Big(1+(\sum_{n=1}^{\infty}\frac{x(rx)^{n-1}}{n!})r\Big)^{-1} \\
& \in E(R[G]^{\tau}, I).
\end{align*}
Hence $Im(\psi) \subset E(R[G]^{\tau}, I)$. For any $ k \geq 1$ and any $x,y \in p^kI$, 
\[
Exp(x)Exp(y) \equiv Exp(x+y) (\text{mod } U(p^kI) \subset p^{2k}U(I) \subset p^{2k}[R[G]^{\tau},I]).
\]
Also, for any $k,l \geq 1$ and any $ x \in p^kI$, $y \in p^lI$, 
\[
Exp(x)Exp(y) \equiv Exp(y)Exp(x) (\text{mod } [p^kI, p^lI] \subset p^{k+l}[R[G]^{\tau},I]).
\]
So for any $ l \geq k \geq 1$, and any $x \in p^k[R[G]^{\tau}, I]$ and $y\in p^l[R[G]^{\tau}, I]$, 
\[
\psi(x) \equiv Exp(x) (\text{mod } p^{2k}[R[G]^{\tau}, I])
\]
\begin{equation}
\label{psiequation}
\psi(x+y) \equiv \psi(x)\psi(y) \equiv \psi(x)Exp(y) (\text{mod } p^{k+l}[R[G]^{\tau}, I]).
\end{equation}
For arbitrary $u \in Exp(p[R[G]^{\tau}, I])$, define a sequence $x_0, x_1, x_2, \ldots$ in $[R[G]^{\tau}, I]$ by setting
\[
x_0= Log(u) \in p[R[G]^{\tau}, I]; \qquad x_{i+1} = x_i + Log(\psi(x_i)^{-1}u),
\]
By (\ref{psiequation}), applied inductively for all $i \geq 0$,
\[
\psi(x_i) \equiv u, \quad x_{i+1} \equiv x_i \quad (\text{mod } p^{2+i}[R[G]^{\tau},I]).
\]
So ${x_i}$ converges and $u = \psi(\lim_{i \rightarrow \infty}   x_i)$. This shows that 
\[
Exp(p[R[G]^{\tau}, I]) \subset Im(\psi) \subset E(R[G]^{\tau}, I).
\]

Now define subgroups $D_k$, for all $ k \geq 0$, by setting
\[
D_k = \langle rx-xr : x \in I, r \in R[G]^{\tau}, rx, xr \in IJ_R^k \rangle \subset [R[G]^{\tau}, I]\cap IJ_R^k.
\]
By the hypothesis on $I$, for all $k \geq 0$, 
\begin{align*}
U(IJ^k_R) &= \sum_{m,n \geq 1} \frac{1}{m+n} [(IJ_R^k)^m, (IJ_R^k)^n] \\
&= \langle [r, \frac{\xi^n}{n}s] : n \geq 2, \xi r, \xi s \in IJ_R^k , \xi^nrs, \xi^nsr \in (IJ_R^k)^n \subset nIJ_R^{k+1} \rangle \\
& \subset D_{k+1}.
\end{align*} 
This shows that $Exp(D_k) \subset Exp([R[G]^{\tau}, I])$ are both normal subgroups of $(R[G]^{\tau})^{\times}$. Also, for any $x, y \in IJ_R^k$,
\begin{equation}
\label{uik}
Exp(x)Exp(y) \equiv Exp(x+y) (\text{mod } Exp(U(IJ_R^k)) \subset Exp(D_{k+1})).
\end{equation}
For any $k \geq 0$ and any $x \in D_k$, if we write $x = \sum (r_ix_i - x_ir_i)$, where $r_i \in R[G]^{\tau}, x_i \in I$ and $r_ix_i, x_ir_i \in IJ_R^k$, then by above 
\begin{align*}
Exp(x) & \equiv \prod (Exp(r_ix_i)Exp(x_ir_i)^{-1}) (\text{mod } Exp(D_{k+1})) \\
& \equiv 1 (\text{mod } E(R[G]^{\tau}, I)).
\end{align*}
In other words, $Exp(D_k) \subset E(R[G]^{\tau}, I)Exp(D_{k+1})$ for all $k \geq 0$. But for a large enough $k$, $D_k \subset p[R[G]^{\tau}, I]$. Hence we get 
\[
Exp(R[G]^{\tau}, I]) \subset Exp(D_0) \subset E(R[G]^{\tau}, I)Exp(p[R[G]^{\tau}, I]) \subset E(R[G]^{\tau}, I).
\]
 \qed

\begin{proposition} For any ideal $I \subset J_R$ of $R[G]^{\tau}$, the $p$-adic logarithm $Log(1+x)$, for any $x\in I$, induces a unique homomorphism 
\[
log_I : K_1(R[G]^{\tau}, I) \rightarrow (I/[R[G]^{\tau},I])\otimes_{\mathbb{Z}_p}\mathbb{Q}_p.
\]
If, furthermore, $I \subset \xi R[G]^{\tau}$ for some central $\xi$ such that $\xi^p \in p\xi R[G]^{\tau}$, then the logarithm induces a homomorphism 
\[
log_I : K_1(R[G]^{\tau}, I) \rightarrow I/[R[G]^{\tau},I],
\]
and $log_I$ is an isomorphism if $I^p \subset pIJ_R$.
\label{logdefn}
\end{proposition}  

\noindent{\bf Proof:} This is an analogue of theorem 2.8 of R. Oliver \cite{Oliver:1988}. By the previous lemma 
\begin{equation}
\label{log}
L : 1+I \xrightarrow{Log} I[\frac{1}{p}] \xrightarrow{proj} I[\frac{1}{p}]/[R[G]^{\tau}, I[\frac{1}{p}]],
\end{equation}
is a homomorphism.  
For each $ n \geq 1$, let
\begin{equation}
\label{trlog}
Tr_n: M_n(I[\frac{1}{p}])/[M_n(R[G]^{\tau}[\frac{1}{p}]), M_n(I[\frac{1}{p}])] \rightarrow I[\frac{1}{p}]/[R[G]^{\tau}[\frac{1}{p}], I[\frac{1}{p}]]
\end{equation}
be the homomorphism induced by the trace map. Then (\ref{log}), applied to the ideal $M_n(I) \subset M_n(R[G]^{\tau})$, induces a homomorphism 
\begin{align*}
L_n: 1+M_n(I) = GL_n(R[G]^{\tau},I) & \xrightarrow{Log} M_n(I[\frac{1}{p}])/[M_n(R[G]^{\tau}[\frac{1}{p}]), M_n(I[\frac{1}{p}])] \\ & \xrightarrow{Tr_n} I[\frac{1}{p}]/[R[G]^{\tau}, I[\frac{1}{p}]].
\end{align*}
For any $n$, and any $u \in 1+M_n(I)$ and $r \in GL_n(R[G]^{\tau})$. 
\[
L_n([r,u]) = L_n(rur^{-1}) - L_n(u) = Tr_n(rLog(u)r^{-1}) - Tr_n(Log(u)) = 0.
\]
So $L_{\infty} = \cup L_n$ factors through a homomorphism 
\begin{align*}
log_I : K_1(R[G]^{\tau}, I) = GL(R[G]^{\tau}, I)/[ & GL(R[G]^{\tau}), GL(R[G]^{\tau}, I)] \\&\rightarrow I[\frac{1}{p}]/[R[G]^{\tau}[\frac{1}{p}], I[\frac{1}{p}]].
\end{align*}

If $I \subset \xi R[G]^{\tau}$, for some central element $\xi$, such that $\xi^p \in p\xi R[G]^{\tau}$, then the same argument by the second part of the above lemma gives a homomorphism 
\[
log_I : K_1(R[G]^{\tau}, I) \rightarrow I/[R[G]^{\tau}, I].
\]
If in addition, $I^p \subset pIJ_R$, then $Log$ is bijective and $Log^{-1}([R[G]^{\tau}, I]) \subset E(R[G]^{\tau}, I])$, by the last part of lemma above. Hence $log_I$ is an isomorphism. 

\qed

\begin{proposition} In the case when $R = \Lambda_O(\Gamma^{p^e})$, the $p$-adic logarithm $Log(1+x)$, for any $x \in J_R$, induces a unique homomorphism 
\[
log : K_1(R[G]^{\tau}) \rightarrow R[Conj(G)]^{\tau}\otimes_{\mathbb{Z}_p} \mathbb{Q}_p.
\]
\end{proposition}

\noindent{\bf Proof:} In this case $R[G]^{\tau}/J_R \cong \mathbb{F}_p$. Hence we have the following exact sequence
\[
1 \rightarrow K_1(R[G]^{\tau}, J_R) \rightarrow K_1(R[G]^{\tau}) \rightarrow K_1(\mathbb{F}_p) = \mathbb{F}_p^{\times}
\]
Hence $log_{J_R}$ extends uniquely to a map 
\begin{equation}
\label{fp}
log : K_1(R[G]^{\tau}) \rightarrow R[Conj(G)]^{\tau} \otimes_{\mathbb{Z}_p}\mathbb{Q}_p. 
\end{equation}

\qed

\begin{remark} The analogue of the equation (\ref{fp}) in the case when $R = \widehat{\Lambda_{O}(\Gamma^{p^e})_{(p)}}$ is
\[
K_2(Q(\mathbb{F}_p[[\Gamma]])) \rightarrow K_1(R[G]^{\tau}, J_R) \rightarrow K_1(R[G]^{\tau}) \rightarrow K_1(Q(\mathbb{F}_p[[\Gamma]])).
\]
Since $K_1(Q(\mathbb{F}_p[[\Gamma]]))$ is not a torsion group, it is not clear if we can extend the homomorphism $log_{J_R}$ to $K_1(R[G]^{\tau})$ is any canonical fashion. However, for our purposes $log_{J_R}$ suffices.
\end{remark}

\subsection{Integral logarithm} Again we follow R. Oliver \cite{Oliver:1988} to construct the integreal logarithm homomorphism from $K_1$ of Iwasawa algebras. This is  a straightforward generalisation of the integral logarithm on $K_1$ of $p$-adic group rings of finite groups contructed by R. Oliver and L. Taylor. 

\noindent{\bf Assumption:} From now on, throughout this paper, we assume that $O$ is a ring of integers in a finite unramified extension of $\mathbb{Q}_p$. \\

\noindent As before $R$ denotes either $\Lambda_O(\Gamma^{p^e})$ or $\widehat{\Lambda_O(\Gamma^{p^e})_{(p)}}$. 

\begin{definition} Let $\varphi$ be the map on $R$ induced by the Frobenius map on $O$ and the $p$-power map on $\Gamma^{p^e}$. We extend this to a map, still denoted by $\varphi$, to $R[Conj(G)]^{\tau}$ by mapping $[\overline{g}]$ to $[\overline{g}^p]$ (remember that $\overline{g}^p$ is not the same as $\overline{g^p}$). 
\end{definition}

\subsubsection{Integral logarithm for $\Lambda_O(\mathcal{G})$} 

\begin{definition} For a finite group $P$ we define the $SK_1(O[P])$ to be the kernel of the map
 \[
 K_1(O[P]) \rightarrow K_1(O[P][\frac{1}{p}]),
 \]
 induced by the natural injection of $O[P]$ into $O[P][\frac{1}{p}]$.
 \end{definition}

For every $n \geq 0$, let $G_n = \mathcal{G}/\Gamma^{p^{e+n}}$. The integral logarithm map defined by R. Oliver and L. Taylor is
\[
L: K_1(O[G_n]) \rightarrow O[Conj(G_n)],
\]
defined as $L = log - \frac{\varphi}{p}log$, where $\varphi: O[Conj(G_n)] \rightarrow O[Conj(G_n)]$ is the map induced by Frobenius on $O$ and the $p$-power map on $G_n$.  The kernel of $L$ is $K_1(O[G_n])_{tor}$ which is equal to $\mu(O) \times G_n^{ab} \times SK_1(O[G_n])$ by a theorem of G. Higman \cite{Higman:1940} and C.T.C. Wall \cite{Wall:1974}. Here $\mu(O)$ denotes the torsion subgroup of $O^{\times}$.
By theorem 7.1 of R. Oliver \cite{Oliver:1988}, $SK_1(O[G_n]) \cong SK_1(O[G])$. 

\begin{definition}  We define the integral logarithm map $L$ on $K_1(\Lambda_O(\mathcal{G}))$ by
 \[
 L(x)= log(x) - \frac{\varphi}{p} log(x) 
 \]
 \end{definition}
 
 \begin{lemma} The image of $L$ is contained in $\Lambda_O(\Gamma^{p^e})[Conj(G)]^{\tau}$. 
\end{lemma}
 
 \noindent{\bf Proof:} For every $n \geq 0$, we have the following commutative diagram
\[
\xymatrix{ K_1(\Lambda_O(\mathcal{G})) \ar[rr]^{L} \ar[d] & & \Lambda_O(\Gamma^{p^e})[Conj(G)]^{\tau}[\frac{1}{p}] \ar[d] \\ 
K_1(O[G_n]) \ar[rr]_{L} & & O[Conj(G_n)][\frac{1}{p}]}
\]
where the vertical arrows are induced by natural projections. By theorem 6.2 of R. Oliver \cite{Oliver:1988}, the image of $L$ is contained in $O[Conj(G_n)]$, for every $n \geq 0$. If an element $x \in \Lambda_O(\Gamma^{p^e})[Conj(G)]^{\tau}[\frac{1}{p}]$ maps into $O[Conj(G_n)]$ for every $ n \geq 0$, then $x$ actually belongs to $\Lambda_O(\Gamma^{p^e})[Conj(G)]^{\tau}$. Hence the lemma. 
\qed

\begin{lemma} The kernel of the map
\[
L : K_1(\Lambda_O(\mathcal{G})) \rightarrow \Lambda_O(\Gamma^{p^e})[Conj(G)]^{\tau},
\]
is $\mu(O) \times \mathcal{G}^{ab} \times SK_1(O[G])$. 
\end{lemma}
\noindent{\bf Proof:} From the proof of previous lemma it is clear that the kernel of $L$ is equal to inverse limit of kernels of $L$ acting on $K_1(O[G_n])$. Hence
\[
Ker(L) = \ilim{n} (\mu(O) \times G_n^{ab} \times SK_1(O[G])) = \mu(O) \times \mathcal{G}^{ab} \times SK_1(O[G]).
\]
\qed
  
\begin{lemma} The cokernel of the map
\[
L : K_1(\Lambda_O(\mathcal{G})) \rightarrow \Lambda_O(\Gamma^{p^e})[Conj(G)]^{\tau},
\]
is $\langle (-1)^{p-1} \rangle \times \mathcal{G}^{ab}$.
\label{integrallogsequence}
\end{lemma}
\noindent{\bf Proof:} We have the following exact sequence by theorem 6.6 of R. Oliver \cite{Oliver:1988}
\[
K_1(O[G_n]) \xrightarrow{L} O[Conj(G_n)] \xrightarrow{\omega} \langle (-1)^{p-1} \rangle \times G_n^{ab} \rightarrow 1.
\]
Here the map $\omega$ is defined by
\[
\sum_{g \in Conj(G_n)} a_gg = \prod_{g} ((-1)^{p-1}g)^{tr_{O/\mathbb{Z}_p}(a_g)}.
\]
Hence cokernel of the map
\[
L : K_1(\Lambda_O(\mathcal{G})) \rightarrow \Lambda_O(\Gamma^{p^e})[Conj(G)]^{\tau},
\]
is $\langle (-1)^{p-1} \rangle \times \mathcal{G}^{ab}$. We denote the map $\Lambda_O(\Gamma^{p^e})[Conj(G)]^{\tau} \rightarrow \langle (-1)^{p-1} \rangle \times \mathcal{G}^{ab}$ by $\omega$.
\qed

\subsubsection{Integral logarithm for $\widehat{\Lambda_O(\mathcal{G})_S}$}

Let $J$ denote the kernel of the natural surjection
\[
\widehat{\Lambda_O(\mathcal{G})_S} \rightarrow \widehat{\Lambda_O(\Gamma)_{(p)}}.
\]
Since $\widehat{\Lambda_O(\mathcal{G})_S}$ is local, the following maps are surjective
\[
\widehat{\Lambda_O(\mathcal{G})_S}^{\times} \rightarrow K_1(\widehat{\Lambda_O(\mathcal{G})_S}),
\]
and 
\[
1+J \rightarrow K_1(\widehat{\Lambda_O(\mathcal{G})_S}, J).
\]
We have an exact sequence 
\[
1 \rightarrow 1+J \rightarrow \widehat{\Lambda_O(\mathcal{G})_S}^{\times} \rightarrow \widehat{\Lambda_O(\Gamma)_{(p)}}^{\times} \rightarrow 1,
\]
which splits using the distinguished embedding $\Gamma \hookrightarrow \mathcal{G}$. Hence any element $x \in \widehat{\Lambda_O(\mathcal{G})_S}^{\times}$ can be written uniquely as $x = uy$, for $u \in 1+J$ and $y \in \widehat{\Lambda_O(\Gamma)_{(p)}}^{\times}$. As a result we get the following
 
 \begin{lemma} Every $x \in K_1(\widehat{\Lambda_O(\mathcal{G})_S})$ can be written as a product $x = uy$, where $ y \in K_1(\widehat{\Lambda_O(\Gamma)_{(p)}})$ and $u \in K_1(\widehat{\Lambda_O(\mathcal{G})_S}, J)$.
\end{lemma}

\begin{lemma} For any $y \in K_1(\widehat{\Lambda_O(\Gamma)_{(p)}})$, we have
\[
\frac{y^p}{\varphi(y)} \equiv 1 (\text{mod } p\widehat{\Lambda_O(\Gamma)_{(p)}})
\]
Hence the logarithm of $\frac{y^p}{\varphi(y)}$ is defined.
\end{lemma}
\noindent{\bf Proof:} The ring $\widehat{\Lambda_O(\Gamma)_{(p)}}/p\widehat{\Lambda_O(\Gamma)_{(p)}}$ is isomorphic to the domain $O/(p)[[\Gamma]]$ (Recall that $(p)$ is the maximal ideal of $O$ since we have assumed that $O$ is unramified extension of $\mathbb{Z}_p$). Let $O/(p) = \mathbb{F}_q$. We also note that $\mathbb{F}_q[[\Gamma]] \cong \mathbb{F}_q[[X]]$. On this ring the $\varphi$ map is the one induced by $X \mapsto (1+X)^p-1 = X^p$ and by $p$-power Frobenius on $\mathbb{F}_q$. Let $ \overline{y} \in \mathbb{F}_q[[X]]$. Write $\overline y = \sum_{i=0}^{\infty} a_iX^i$. Then
\[
\overline{y}^p = (\sum_{i} a_iX^i)^p = \sum_i a_i^p X^{ip} = \varphi(y),
\]
Hence the lemma.
\qed
 
\begin{definition} We define the integral logarithm $L$ on $K_1(\widehat{\Lambda_O(\mathcal{G})_S})$ as follows: write any $x \in K_1(\widehat{\Lambda_O(\mathcal{G})_S})$ as $x = uy$, with $y \in K_1(\widehat{\Lambda_O(\Gamma)_{(p)}})$ and $u \in K_1(\widehat{\Lambda_O(\mathcal{G})_S}, J)$. Define
\[
L(x) = L(uy) = L(u)+L(y) = log(u)- \frac{\varphi}{p}log(u) + \frac{1}{p}log(\frac{y^p}{\varphi(y)}).
\]
\end{definition}

\begin{proposition} $L$ induces a homomorphism 
\[
L: K_1(\widehat{\Lambda_O(\mathcal{G})_S}) \rightarrow \widehat{\Lambda_O(\Gamma^{p^e})_{(p)}}[Conj(G)]^{\tau}.
\]
\end{proposition}
\noindent{\bf Proof:} Let $x \in K_1(\widehat{\Lambda_O(\mathcal{G})_S})$. Write it as $x=uy$, with  $y \in K_1(\widehat{\Lambda_O(\Gamma)_{(p)}})$ and $u \in K_1(\widehat{\Lambda_O(\mathcal{G})_S}, J)$. We will show that $L(u)$ and $L(y)$ lie in $\widehat{\Lambda_O(\Gamma^{p^e})_{(p)}}[Conj(G)]^{\tau}$. The fact that $L(y) = \frac{1}{p}log(\frac{y^p}{\varphi(y)})$ lies in $\widehat{\Lambda_O(\Gamma^{p^e})_{(p)}}$ follows from proposition \ref{logdefn}. Now let $u=1-v$ and consider $L(1-v)$
\begin{align*}
L(1-v) &= -(v+\frac{v^2}{2} +\frac{v^3}{3}+ \cdots)+(\frac{\varphi(v)}{p}+\frac{\varphi(v)}{2p}+\cdots) \\
& \equiv \sum_{k=1}^{\infty} \frac{1}{pk}(v^{pk}-\varphi(v^k)) \quad (\text{mod } \widehat{\Lambda_O(\Gamma^{p^e})_{(p)}}[Cong(G)]^{\tau})
\end{align*}
It suffices that $pk| (v^{pk} - \varphi(v^k))$ for all $k$; which will follow from 
\[
p^n|(v^{p^n}-\varphi(v^{p^{n-1}})),
\]
for all $n \geq 1$. Write $v = \sum_{g \in G} r_g\bar{g}$, where $r_g \in \widehat{\Lambda_O(\Gamma^{p^e})_{(p)}}$. Set $ q =p^n$ and consider a typical term in $v^q$:
\[
r_{g_1}\cdots r_{g_q} \bar{g}_1\cdots \bar{g}_q.
\]
Let $\mathbb{Z}/q\mathbb{Z}$ acts cyclically permuting the $g_i$'s, so that we get a total of $p^{n-t}$ conjugate terms, where $p^t$ is the number of cyclic permutations leaving each term invariant. Then $\bar{g}_1\cdots \bar{g}_q$ is a $p^t$-th power, and the sum of the conjugate terms has the form
\[
p^{n-t} \hat{r}^{p^t} \hat{g}^{p^t} \in \widehat{\Lambda_O(\Gamma^{p^e})_{(p)}}[Conj(G)]^{\tau},
\]
where $\hat{r} = \prod_{j=1}^{p^{n-t}} r_{g_j}$ and $\hat{g} = \prod_{j=1}^{p^{n-t}} g_j$. Here we may do the multiplication after rearrangement because the twisting map $\tau$ is symmetric. If $t=0$, then this is a multiple of $p^n$. If $t >0$, then there is a corresponding term $p^{n-t} \hat{r}^{p^{t-1}} \hat{g}^{p^{t-1}}$ in the expansion of $v^{p^{n-1}}$. It remains only to show that
\[
p^{n-t}\hat{r}^{p^t}\hat{g}^{p^t} \equiv p^{n-t} \varphi(\hat{r}^{p^{t-1}}\hat{g}^{p^{t-1}}) = p^{n-t} \varphi(\hat{r}^{p^{t-1}}) \hat{g}^{p^t} (\text{mod } p^n).
\]
But $p^t |(\hat{r}^{p^t} - \varphi(\hat{r}^{p^{t-1}}))$, since $p|(\hat{r}^p - \varphi(\hat{r}))$. 
\qed

\begin{remark} It would be interesting to find the kernel and cokernel of $L$ in this case. Since we do not need it we will not investigate it here.
\end{remark}

\section{The main result}
\subsection{The maps $\theta^G$ and $\theta^G_S$}

 Recall that for every $P \in C(G)$, we denote by $U_P$ the inverse image of $P$ in $\mathcal{G}$. Then there is a map $\theta^G_P$ given by the norm
\[
\theta^G_P: K_1(\Lambda_O(\mathcal{G})) \xrightarrow{norm} K_1(\Lambda_O(U_P)) = \Lambda_O(U_P)^{\times}. 
\]
Let $\theta^G$ be the map
\[
\theta^G = (\theta^G_P)_{P \in C(G)} : K_1(\Lambda_O(\mathcal{G})) \rightarrow \prod_{P \in C(G)} \Lambda_O(U_P)^{\times}.
\]

Similarly, we have maps $\theta^G_S$ and $\widehat{\theta}^G$
\[
\theta^G_S : K_1(\Lambda_O(\mathcal{G})_S) \rightarrow \prod_{P \in C(G)} \Lambda_O(U_P)_S^{\times}.
\]
\[
\widehat{\theta}^G : K_1(\widehat{\Lambda_O(\mathcal{G})_S}) \rightarrow \prod_{P \in C(G)} \widehat{\Lambda_O(U_P)_S}^{\times}. 
\]
Our main theorem describes kernel and image of the map $\theta^G$ and shows that image of $\theta^G_S$ intersected with $\prod_{P \in C(G)} \Lambda_O(U_P)^{\times}$ is exactly the image of $\theta^G$.

\subsection{Relation between the maps $\theta$ and $\beta$} 

\begin{lemma} For any $P \in C(G)$, we have the following diagram 
\[
\xymatrix{ K_1(\Lambda_O(\mathcal{G})) \ar[rr]^{log} \ar[d]_{\theta^G_P} & & \Lambda_O(\Gamma^{p^e})[Conj(G)]^{\tau}[\frac{1}{p}] \ar[d]^{t^G_P} \\
\Lambda_O(U_P)^{\times} \ar[rr]_{log} & & \Lambda_O(\Gamma^{p^e})[P]^{\tau}[\frac{1}{p}]}
\]
In the case of $p$-adic completions the following diagram commutes
\[
\xymatrix{ K_1(\widehat{\Lambda_O(\mathcal{G})_S}, J) \ar[rr]^{log} \ar[d]_{\widehat{\theta^G_P}} & & \widehat{\Lambda_O(\Gamma^{p^e})_{(p)}}[Conj(G)]^{\tau}[\frac{1}{p}] \ar[d]^{t^G_P} \\
K_1(\widehat{\Lambda_O(U_P)_S},J) \ar[rr]_{log} & & \widehat{\Lambda_O(\Gamma^{p^e})_{(p)}}[P]^{\tau}[\frac{1}{p}]}
\]
\end{lemma}
\noindent{\bf Proof:} This is an analogue of theorem 6.8 in R. Oliver \cite{Oliver:1988} which we refer for details. For any $u \in K_1(\Lambda_O(\mathcal{G}))$ (resp. $u \in K_1(\widehat{\Lambda_O(\mathcal{G})_S}, J)$), the result follows from the expressions
\[
log(u) = \lim_{n \rightarrow \infty} \frac{1}{p^n}(u^{p^n}-1); 
\]
and  
\[
\theta^G_P(u) = \lim_{n \rightarrow \infty}(1+t^G_P(u^{p^n}-1))^{1/p^n}.
\]
\[
(\text{resp.  } \widehat{\theta^G_P}(u) = \lim_{n \rightarrow \infty}(1+t^G_P(u^{p^n}-1))^{1/p^n} ).
\]
\qed
 
\begin{definition} For any $P \in C(G)$ with $P \neq \{1\}$, we choose and fix a non-trivial character $\omega$ of $P$ of order $p$. It induces a map on $\Lambda_O(U_P)^{\times}$ and on $\widehat{\Lambda_O(U_P)_S}^{\times}$ given by mapping $g \in U_P$ to $\omega(g)g$. We denote this map by $\omega$ again. We define a map $\alpha_P$ from $\Lambda_O(U_P)^{\times}$ to itself or from $\widehat{\Lambda_O(U_P)_S}^{\times}$ to itself by 
\[
\alpha_P(x) = \frac{x^p}{\prod_{k=0}^{p-1}\omega^k(x)}.
\]
We define $\alpha_{\{1\}}$ by $\alpha_{\{1\}}(x) = x^p/\varphi(x)$. We put $\alpha = (\alpha_P)_{P \in C(G)}$.
\end{definition}

\begin{lemma} We have the following commutative diagram for all $P \in C(G)$ and $P \neq \{1\}$
\[
\xymatrix{ \Lambda_O(U_P)^{\times} \ar[r]^{log} \ar[d]_{\alpha_P} &  \Lambda_O(U_P)[\frac{1}{p}] \ar[d]^{p\eta_P} \\ 
 \Lambda_O(U_P)^{\times} \ar[r]_{log} &  \Lambda_O(U_P)[\frac{1}{p}]}.
\]
And in the case of $p$-adic completion the following diagram commutes
\[
\xymatrix{  K_1(\widehat{\Lambda_O(U_P)_S}, J) \ar[rr]^{log} \ar[d]_{\alpha_P} & &  \widehat{\Lambda_O(U_P)_S}[\frac{1}{p}] \ar[d]^{p\eta_P}\\
K_1(\widehat{\Lambda_O(U_P)_S},J) \ar[rr]_{log} & &  \widehat{\Lambda_O(U_P)_S}[\frac{1}{p}] }
\]

\end{lemma}

\noindent{\bf Proof:} Note that the following diagram commutes
\[
\xymatrix{ \Lambda_O(U_P)^{\times} \ar[r]^{log} \ar[d]_{\omega} & \Lambda_O(U_P)[\frac{1}{p}] \ar[d]^{\omega} \\
\Lambda_O(U_P)^{\times} \ar[r]_{log} & \Lambda_O(U_P)[\frac{1}{p}] }
\]
and
\[
\eta_P = \frac{1}{p}(p - \sum_{k=0}^{p-1} \omega^k).
\]
Hence we get the first claim of the lemma. The second one is similar. 
\qed

\begin{definition} We define the map $v^G_P : \prod_{C\in C(G)} \Lambda_O(U_C) \rightarrow  \Lambda_O(U_P)$ is given by 
\[
v^G_P((x_C)) =p(\sum ver^{P'}_P(x_{P'}))
\]
where the sum ranges over all $P' \in C(G)$, such that $P'^p=P$ and $P \neq P'$. The empty sum is taken to be 0. Here $ver^{P'}_P$ is the map induced by the transfer homomorphism (equal to the $p$-power map) from $P'$ to $P$. Put $v^G = (v^G_P)_P$. We denote the analogous map in the case of $p$-adic completions again by $v^G$.
\end{definition}

\begin{lemma} For $P \neq \{1\}$, the following diagram commutes 
\[
\xymatrix{ \Lambda_O(\Gamma^{p^e})[Conj(G)]^{\tau} \ar[rr]^{\varphi} \ar[d]_{\beta^G} & & \Lambda_O(\Gamma^{p^e})[Conj(G)]^{\tau} \ar[d]^{\beta^G_P} \\
\prod_{C\in C(G)} \Lambda_O(U_C) \ar[rr]_{v^G_P} & &  \Lambda_O(U_P)}
\]
In the case $P = \{1\}$, the following diagram commutes
\[
\xymatrix{ \Lambda_O(\Gamma^{p^e})[Conj(G)]^{\tau} \ar[rr]^{\varphi} \ar[d]_{\beta^G} & & \Lambda_O(\Gamma^{p^e})[Conj(G)]^{\tau} \ar[d]^{\beta^G_P} \\
\prod_{C\in C(G)} \Lambda_O(U_C) \ar[rr]_{\varphi + v^G_P} & &  \Lambda_O(U_P)}
\]
Here $\varphi$ in the lower row is just the map $\varphi: \Lambda_O(U_P) \rightarrow \Lambda_O(U_P)$. And the analogous result holds for the $p$-adic completions. 
\end{lemma} 

\noindent{\bf Proof:} We assume that $P \neq \{1\}$. For $g \in G$, we must show that $\beta^G_P(\varphi([\overline{g}])) = v^G_P(\beta^G([\overline{g}]))$. We start with the left hand side 
\[
\beta^G_P(\varphi([\overline{g}])) = \eta_P(t^G_P([\overline{g}^p])),
\]
which is non-zero if and only if some conjugate of $g^p$ generates $P$. Hence we may assume that $g^p$ itself is a generator of $P$. Then
\begin{align*}
v^G_P(\beta^G([\overline{g}])) & = p(\sum_{P' \text{ s.t. } P'^p=P} ver^{P'}_P(\beta^G_{P'}([\overline{g}]))) \\
& = p(\sum_{P' \text{ s.t. } P'^p=P} ver^{P'}_P(t^G_{P'}([\overline{g}]))) \\
& = p(\sum_{P' \text{ s.t. } P'^p=P} ver^{P'}_P(\sum_{x \in C(G,P')} \{ [\overline{x^{-1}gx}] | x^{-1}gx \in P'\})) \\
&= p(\sum_{P' \text{ s.t. } P'^p=P} \sum_{x \in C(G,P')} \{[\overline{x^{-1}g^px}] | x^{-1}gx \in P'\}) \\
&= \sum_{P' \text{ s.t. } P'^p=P}  (\sum_{x \in C(G,P)} \{[\overline{x^{-1}g^px}] | x^{-1}gx \in P'\}) \\
& = \sum_{x \in C(G,P)} \sum_{P' \text{ s.t. } P'^p=P} \{[\overline{x^{-1}g^px}] | x^{-1}gx \in P'\} \\
& = \sum_{x \in C(G,P)} \{[\overline{x^{-1}g^px}] | x^{-1}g^px \in P\} \\
& = \beta^G_P(\varphi([\overline{g}])).
\end{align*}

\qed

\begin{definition} We define the map $u^G_P: \prod_{C \in C(G)} \Lambda_O(U_C)^{\times} \rightarrow \Lambda_O(U_P)^{\times}$ by 
\[
u^G_P((x_C)) = \prod ver^{P'}_P(x_{P'}),
\]
where the product is taken over all $P' \in C(G)$ such that $P'^p = P$ and $P' \neq P$. The empty product is taken to be 1. Put $u^G = (u^G_P)_P$. We denote the analogous map in the case of $p$-adic completions again by $u^G$.
\end{definition}

\begin{lemma} The following diagram commutes
\[
\xymatrix{ \prod_{P \in C(G)} \Lambda_O(U_P)^{\times} \ar[r]^{log} \ar[d]_{u^G} & \prod_{P \in C(G)} \Lambda_O(U_P)[\frac{1}{p}] \ar[d]^{\frac{1}{p}v^G} \\
\prod_{P \in C(G)} \Lambda_O(U_P)^{\times} \ar[r]_{log} & \prod_{P \in C(G)} \Lambda_O(U_P)[\frac{1}{p}]}
\]
And similarly for the $p$-adic completions.
\end{lemma}

\noindent{\bf Proof:} This is because $log$ is natural with respect to the homomorphisms of Iwasawa algebras induced by group homomorphisms. 
\qed

\begin{proposition} The additive and the multiplicative sides are related by the following formula:  
\[
\beta^G_P(L(x)) = \frac{1}{p}log(\frac{\alpha_P(\theta^G_P(x))}{u^G_P(\alpha(\theta^G(x)))})
\]
Analogous relation holds in the case of $p$-adic completions.
\label{relation}
\end{proposition}
\noindent{\bf Proof:} First we assume $P \neq \{1\}$. Consider the left hand side of the equation
\begin{align*} 
\beta^G_P(L(x)) &= \eta_P(t^G_P(log(x) - \frac{\varphi}{p} log(x))) \\
& =  \eta_P(log(\theta^G_P(x))) - \beta^G_P(\frac{\varphi}{p}log(x)) \\
& = \frac{1}{p}log(\alpha_P(\theta^G_P(x))) - \frac{1}{p}v^G_P(\beta^G(log(x))) \\
& = \frac{1}{p} log(\alpha_P(\theta^G_P(x)) - \frac{1}{p}log(u^G_P(\alpha(\theta^G(x)))) \\
& = \frac{1}{p} log\Big( \frac{\alpha_P(\theta^G_P(x))}{u^G_P(\alpha(\theta^G(x)))} \Big).
\end{align*}

Now assume that $P = \{1\}$. In this case $\beta^G_P = t^G_P$ and
\begin{align*}
\beta^G_P(L(x)) &= t^G_P(log(x)) - \frac{1}{p}t^G_P(\varphi(log(x))) \\
&= log(\theta^G_P(x)) - \frac{1}{p}(\varphi+v^G_P)(\beta^G(log(x))) \\
&= log(\theta^G_P(x)) - \frac{1}{p^2}(p\varphi+v^G_P)(log(\alpha(\theta^G(x)))) \\
& = \frac{1}{p} log\Big( \frac{\theta^G_P(x)^p}{\varphi(\theta^G_P(x))u^G_P(\alpha(\theta^G(x)))}\Big) \\
& = \frac{1}{p} log\Big(\frac{\alpha_P(\theta^G_P(x))}{u^G_P(\alpha(\theta^G(x)))}\Big)
\end{align*}

\qed
 
\subsection{Statement of the main theorems}

\begin{definition} Let $\Phi^G$ (resp. $\Phi^G_S$ and $\widehat{\Phi^G_S}$) be the subgroup of $\prod_{P\in C(G)}\Lambda_O(U_P)^{\times}$ (resp. $\prod_{P\in C(G)} \Lambda_O(U_P)_S^{\times}$ and $\prod_{P \in C(G)}\widehat{\Lambda_O(U_P)_S}^{\times}$) consisting of tuples $(x_P)$ satisfying \\
M1. For any $P \leq P_1$ in $C(G)$, we have 
\[
nr^{P_1}_P(x_{P_1}) = x_P.
\]
M2. $(x_P)$ is fixed under conjugation action by every $g \in G$. \\
M3. $(x_P)$ satisfies the following congruence
\[
\alpha_P(x_P) \equiv   u^G_P(\alpha(x_{P'})) (\text{mod } pT_P) (\text{ resp. } pT_{P,S} \text{ and } p\widehat{T_P}),
\]

\end{definition} 
 
\begin{theorem}The map $\theta^G$ induces an isomorphism 
\[
K_1(\Lambda_O(\mathcal{G}))/SK_1(\Lambda_O(\mathcal{G})) \cong \Phi^G.
\]
\label{theorem1}
\end{theorem}

\begin{theorem} The image of $\theta^G_S$ is contained in $\Phi^G_S$. Moreover, 
\[
\Phi^G_S \cap \prod_{P \in C(G)} \Lambda_O(U_P)^{\times} = Im(\theta^G).
\]
\label{theorem2}
\end{theorem}

\subsection{The proof} In this section we again let $R$ denote either the ring $\Lambda_O(\Gamma^{p^e})$ or the ring $\widehat{\Lambda_O(\Gamma^{p^e})_{(p)}}$. 

\begin{lemma} Let $x \in K_1(R[G]^{\tau})$. Then for every $P \in C(G)$, we have
\[
\alpha_P(\theta^G_P(x)) \equiv u^G_P(\alpha(\theta^G(x))) (\text{mod } p)
\]
In particular, $log$ of $\frac{\alpha_P(\theta^G_P(x))}{u^G_P(\alpha(\theta^G(x)))}$ is defined.
\label{weakcong}
\end{lemma}

\noindent{\bf Proof:} If $P$ is not the trivial group, then $\alpha_P(\theta^G_P(x)) \equiv 1 (\text{mod }p)$ and $u^G_P(\alpha(\theta^G(x))) \equiv 1 (\text{mod }p)$, hence we get the required congruence. Now assume that $P =\{1\}$. In this case
\begin{align*}
u^G_{\{1\}}(\alpha(\theta^G(x))) & = \Big(\prod_{P' \text{ s.t. } P' \neq \{1\}, P'^p= \{1\}}ver^{P'}_P(\alpha_{P'}(\theta^G_{P'}(x)))\Big) \\
& \equiv 1 (\text{ mod } p)
\end{align*}
And 
\[
\theta^G_{\{1\}}(x)^p \equiv \varphi(\theta^G_{\{1\}}(x)) (\text{mod } p).
\]
i.e.
\[
\alpha_{\{1\}}(\theta^G_{\{1\}}(x)) \equiv 1 (\text{mod } p).
\]
\qed

\begin{lemma} The image of $\theta^G$ is contained in $\Phi^G$. Similarly, image of $\widehat{\theta^G_S}$ is contained in $\widehat{\Phi^G_S}$. 
\label{multiplicativecontain}
\end{lemma}

\noindent{\bf Proof:} Here we write a proof for the first claim of the lemma only. The proof for the second claim is exactly the same. We must show that for any $ x \in K_1(\Lambda_O(\mathcal{G})),$ the element $\theta^G(x)$ satisfies M1, M2 and M3. \\

\noindent M1 is clear because the following digram commutes
\[
\xymatrix{ K_1(\Lambda_O(\mathcal{G})) \ar[r]^{\theta^G_{P_1}} \ar[rd]_{\theta^G_P} & \Lambda_O(U_{P_1})^{\times} \ar[d]^{\theta^{P_1}_P}\\
& \Lambda_O(U_P)^{\times}}
\]

\noindent M2. Pick a lift $\tilde{x}$ of $x$ in $\Lambda_O(\mathcal{G})^{\times}$. Recall how the norm map $\theta^G_P$ is defined for any $P \in C(G)$. The ring $\Lambda_O(\mathcal{G})$ is a free $\Lambda_O(U_P)$-module of rank $[G:P]$. We take the $\Lambda_O(U_P)$-linear map on $\Lambda_O(\mathcal{G})$ induced by multiplication by $x$ on the right. Let $A_P(x)$ be the matrix (with entires in $\Lambda_O(U_P)$) of this map with respect to some basis $C(G,P)$. For any $g \in G$, the set $g C(G,P) g^{-1}$ is a basis for the $\Lambda_O(U_{gPg^{-1}})$-module $\Lambda_O(\mathcal{G})$. The matrix for the $\Lambda_O(U_{gPg^{-1}})$-linear map on $\Lambda_O(\mathcal{G})$ induced by multiplication on the right by $x$ with respect to this basis is $gA_P(x)g^{-1}$. Hence
\[
g\theta^G_P(x)g^{-1} = \theta^G_{gPg^{-1}}(x).
\]

\noindent M3. By lemma \ref{weakcong} and proposition \ref{relation} we have
\[
\beta^G_P(L(x)) = \frac{1}{p}log\Big( \frac{\alpha_P(\theta^G_P(x))}{u^G_P(\alpha(\theta^G(x)))} \Big).
\]
By lemma \ref{additivecontain} we get that $p\beta^G_P(L(x)) \in pT_P$. But $log$ induces an isomorphism between $1+pT_P$ and $pT_P$. Hence M3 is satisfied. 

\qed

\begin{lemma} There is a map $\mu(O) \times \mathcal{G}^{ab} \rightarrow \prod_{P\in C(G)} \Lambda_O(U_P)^{\times}$ which is injective and fits in the following commutative diagram
\[
\xymatrix{ \mu(O)\times \mathcal{G}^{ab} \ar@{^{(}->}[r] \ar@{=}[d] & K'_1(\Lambda_O(\mathcal(G))) \ar[d]_{\theta^G} \\
\mu(O)\times \mathcal{G}^{ab} \ar@{^{(}->}[r] & \prod_{P\in C(G)} \Lambda_O(U_P)^{\times} }
\]
where $K'_1(\Lambda_O(\mathcal{G}))$ denotes the quotient $K_1(\Lambda_O(\mathcal{G}))/SK_1(O[G])$.
\end{lemma}

\noindent{\bf Proof:} We define the claimed map as the composition
\[
\mu(O) \times \mathcal{G}^{ab} \xrightarrow{id \times ver} \mu(O) \times U_P \hookrightarrow \Lambda_O(U_P)^{\times}.
\]
To prove the injectivity it suffices to show that
\[
\mathcal{G}^{ab} \rightarrow \prod_{P\in C(G)} U_P
\]
is injective. We use induction on the order of $H$. Take a central commutator $z$ of $\mathcal{G}$ of order $p$ such that $z \in H$. Let $\overline{\mathcal{G}} = \mathcal{G}/\langle z \rangle$. We use the notation $\overline{U_P}$ for image of $U_P$ in $\overline{\mathcal{G}}$. Consider the commuting diagram
\[
\xymatrix{ 1 \ar[r] \ar[d] & \mathcal{G}^{ab} \ar@{=}[r] \ar[d] & \overline{\mathcal{G}}^{ab} \ar[d] \\
K \ar[r] & \prod_{P\in C(G)} U_P \ar[r] & \prod_{P\in C(G)} \overline{U_P}}
\]
where $K$ is by definition the kernel of the map $\prod_{P\in C(G)} U_P \rightarrow \prod_{P\in C(G)} \overline{U_P}$. Since the first and the last vertical maps are injective, the middle one must also be so. Since the norm homomorphism on $K_1$ groups is same as the transfer when restricted to $\mathcal{G}^{ab}$, we also get the commutativity of the diagram in the lemma.

\qed

\begin{lemma} We denote $\psi^G_{\Lambda_O(\Gamma^{p^e})}$ from the additive theorem by $\psi^G$. Then there is a surjective map 
\[
\omega : \psi^G \rightarrow \langle (-1)^{p-1} \rangle \times \mathcal{G}^{ab},
\]
which makes the following diagram commute
\[
\xymatrix{ \Lambda_O(\Gamma^{p^e})[Conj(G)]^{\tau} \ar[r]^{\omega} \ar[d]_{\beta^G} & \langle (-1)^{p-1} \rangle \times \mathcal{G}^{ab} \ar[d] \\
\psi^G \ar[r]_{\omega} & \langle (-1)^{p-1} \rangle \times \mathcal{G}^{ab} }
\]
where the map $\omega$ in the top row is the one from the proof of lemma \ref{integrallogsequence}.
\end{lemma}

\noindent{\bf Proof:} This is a trivial corollary of theorem \ref{additivetheorem}. 

\qed

\begin{lemma} The map $\mathcal{L} : \Phi^G \rightarrow \psi^G$, defined by 
\[
\mathcal{L}_P((x_C)) = \frac{1}{p}log\Big(\frac{\alpha_P(x_P)}{u^G_P(\alpha((x_C)))}\Big)
\]
sits in the following exact sequence
\[
1 \rightarrow \mu(O) \times \mathcal{G}^{ab} \rightarrow \Phi^G \xrightarrow{\mathcal{L}} \psi^G \xrightarrow{\omega} \langle (-1)^{p-1} \rangle \times \mathcal{G}^{ab} \rightarrow 1.
\]
\end{lemma}

\noindent{\bf Proof:} It is clear that the map $\mathcal{L}$ is well-defined and its image is contained in $\psi^G$. To show that it sits in the claimed exact sequence, we use induction on the order of $H$. The only non-trivial part is the exactness at $\psi^G$. Choose a central commutator $z \in H$ of order $p$. Let $\overline{\mathcal{G}} = \mathcal{G}/\langle z \rangle$. Then we have the following commutative diagram
\[
\xymatrix{ \mu(O) \times \mathcal{G}^{ab} \ar[r] \ar@{=}[d] & \Phi^G \ar[r]^{\mathcal{L}} \ar[d]_{\pi} & \psi^G \ar[rr]^{\omega} \ar[d]_{\pi} & & \langle (-1)^{p-1} \rangle \times \mathcal{G}^{ab} \ar@{=}[d] \\
\mu(O) \times \overline{\mathcal{G}}^{ab} \ar[r] & \Phi^{\overline{G}} \ar[r]_{\mathcal{L}} & \psi^{\overline{G}} \ar[rr]_{\omega} & & \langle (-1)^{p-1} \rangle \times \overline{\mathcal{G}}^{ab} }
\]
For any $a = \mathcal{L}(x) \in \psi^G$, we have $\omega(a) = \omega(\pi(a)) = \omega (\mathcal{L}(\pi(x))) = 1$. On the other hand, let $a \in \psi^G$ be in the kernel of $\omega$. Then by the previous lemma, the additive theorem and lemma \ref{integrallogsequence} we get $x \in K_1(\Lambda_O(\mathcal{G}))$ such that $\beta^G(L(x)) = a$. We finish the proof by using the commutativity of the following diagram
\[
\xymatrix{ K_1(\Lambda_O(\mathcal{G}) \ar[rr]^{L} \ar[d]_{\theta^G} & & \Lambda_O(\Gamma^{p^e})[Conj(G)]^{\tau} \ar[d]^{\beta^G} \\
\Phi^G \ar[rr]_{\mathcal{L}} & & \psi^G}
\]
Hence we have $\mathcal{L} (\theta^G(x)) = a$ and we get the exactness at $\psi^G$.	

\qed

\noindent{\bf Proof of theorem \ref{theorem1}:} Recall our notation 
\[
K'_1(\Lambda_O(\mathcal{G})) = K_1(\Lambda_O(\mathcal{G}))/SK_1(O[G]). 
\]
From the previous lemmas we have the following commutative diagram
\[
\xymatrix{ \mu(O) \times \mathcal{G}^{ab} \ar[r] \ar@{=}[d] & K'_1(\Lambda_O(\mathcal{G})) \ar[r] \ar[d]_{\theta^G} & \Lambda_O(\Gamma^{p^e})[Conj(G)]^{\tau} \ar[r] \ar[d]^{\beta^G} & \langle (-1)^{p-1} \rangle \times \mathcal{G}^{ab} \ar@{=}[d]  \\
\mu(O) \times \mathcal{G}^{ab} \ar[r] & \Phi^G \ar[r] & \psi^{G} \ar[r] & \langle (-1)^{p-1} \rangle \times \mathcal{G}^{ab} }
\]
The theorem now follows by the five lemma. 

\qed

\noindent{\bf Proof of theorem \ref{theorem2}:} It follows from lemma \ref{multiplicativecontain} and the fact  
\[
\widehat{\Phi^G_S} \cap \prod_{P \in C(G)} \Lambda_O(U_P)_S^{\times} = \Phi^G_S,
\]
that the image of $\theta^G_S$ is contained in $\Phi^G_S$. Moreover, from 
\[
\Phi^G_S \cap \prod_{P \in C(G)} \Lambda_O(U_P)^{\times} = \Phi^G,
\] 
and theorem \ref{theorem1}, we obtain that 
\[
\Phi^G_S \cap \prod_{P \in C(G)} \Lambda_O(U_P)^{\times} = Im(\theta^G).
\]

\qed

\bibliographystyle{plain}
\bibliography{mybib2}

\end{document}